\documentclass[11pt,sloppy]{article}
 \usepackage{wseas}
  \usepackage{times}
\usepackage{graphicx}
\usepackage{latexsym}
\usepackage{amsfonts}
\usepackage{amsthm}
\usepackage{amsmath}
\usepackage{amssymb}
\usepackage{graphicx}
\usepackage{latexsym}
\usepackage{amsfonts}
\usepackage{amsthm}
\usepackage{amsmath}
\usepackage{amssymb}
\usepackage{multicol}

\newtheorem{thm}{Theorem}[section]
\newtheorem{prop}{Proposition}[section]
\newtheorem{lemma}{Lemma}[section]

\newtheorem{corollary}{Corollary}[thm]

\newcommand{\reals}{\mathbb{R}}

\newcommand{\naturals}{\mathbb{N}}

   \newcommand{\comment}[1]{}

\makeatother

\begin{document}

\title{{\fontsize{16}{20} \bf Magic graphs 
 and the faces of the Birkhoff
polytope}}

\author{{\large MAYA MOHSIN AHMED}  \\  {\large Department of Mathematics} \\
{\large University of California, Davis, CA 95616} \\
{\large UNITED STATES OF AMERICA} \\
{\large e-mail: maya@math.ucdavis.edu \hspace{0.2in} 
 http://www.math.ucdavis.edu/ $\tilde{ }$ maya} 
}

\date{}

\twocolumn[
\begin{@twocolumnfalse} 

\maketitle

\noindent {\bf\normalsize ABSTRACT}\newline Magic labelings of graphs
are studied in great detail by Stanley in \cite{stanley3} and
\cite{stanley4}, and Stewart in \cite{stewart} and
\cite{stewartcomplete}. In this article, we construct and enumerate
magic labelings of graphs using Hilbert bases of polyhedral cones and
Ehrhart quasi-polynomials of polytopes.  We define polytopes of magic
labelings of graphs and digraphs. We give a description of the faces
of the Birkhoff polytope as polytopes of magic labelings of
digraphs. \newline

\noindent {\bf\normalsize KEY WORDS}\newline
Magic graphs; Polyhedral cones; Birkhoff polytope.

\end{@twocolumnfalse}
] 

\section{\Large{Introduction.}} \label{introductionsection}
\begin{figure}
 \begin{center}
     \includegraphics[scale=0.5]{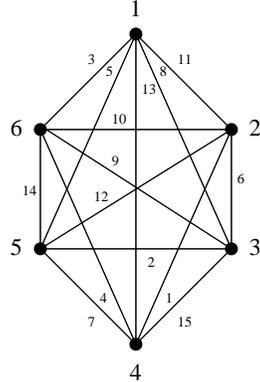}
\caption{A magic labeling of $K_6$ \cite{stewartcomplete}.}  \label{mleg}
 \end{center}
 \end{figure}

Let $G$ be a finite graph. A {\em labeling} of $G$ is an assignment of
a nonnegative integer to each edge of $G$. A {\em magic labeling of
magic sum} $r$ of $G$ is a labeling such that for each vertex $v$ of
$G$ the sum of the labels of all edges incident to $v$ is the magic
sum $r$ (loops are counted as incident only once) \cite{stanley3,
stanley4, stewart, stewartcomplete}. Graphs with a magic labeling are
also called {\em magic graphs} (see Figure \ref{mleg} for an example
of a magic labeling of the complete graph $K_6$ of magic sum $40$).

We define a {\em magic labeling} of a digraph $D$ of {\em magic sum}
$r$ to be an assignment of a nonnegative integer to each edge of $D$,
such that for each vertex $v_i$ of $D$, the sum of the labels of all
edges with $v_i$ as the initial vertex is $r$, and the sum of the
labels of all edges with $v_i$ as the terminal vertex is also $r$.
Thus magic labelings of a digraph is a network flow, where the flow
into and out of every vertex, is the magic sum of the labeling. If we
consider the labels of the edges of $G$ as variables, the defining
magic sum conditions are simply linear equations, and the set of magic
labelings of $G$ becomes the set of integral points inside a {\em
pointed polyhedral cone} $C_G$ \cite{schrijver}. Henceforth, we call
$C_G$ the cone of magic labelings of $G$. A {\em Hilbert basis} of
$C_G$ has the property that any magic labeling of $G$ can be expressed
as a linear combination with nonnegative integer coefficients of the
elements of the Hilbert basis (see \cite{adh}, \cite{ma}, or
\cite{schrijver}).  An {\em irreducible magic labeling} of a graph is
a magic labeling that cannot be written as a sum of two other magic
labelings.  The {\em minimal} Hilbert basis of $C_G$ is the set of all
irreducible magic labelings and is unique
\cite{schrijver}. Henceforth, when we refer to the Hilbert basis in
this article, we mean the minimal Hilbert basis.

Let $v_1, v_2, \dots, v_n$ denote the vertices of $G$ and let
$e_{i_1}, e_{i_2}, \dots, e_{i_{m_i}}$ denote the edges of $G$ that
are incident to the vertex $v_i$ of $G$.  Consider the polytope
\[
{\cal P}_G = \{ L \in C_G \subseteq \reals^{q}, \hspace{0.05 in}
\sum_{j=1}^{m_i}   L(e_{i_j}) = 1 ;  i=1, \dots, n
\}.
\]

We will refer to ${\cal P}_G$ as the polytope of magic labelings of
$G$ and denote $H_G(r)$ to be the number of magic labelings of $G$ of
magic sum $r$.  Then, $H_G(r)$ is the {\em Ehrhart quasi-polynomial}
of ${\cal P}_G$ \cite{stanley}.  

We define a polytope ${\cal P}_D$ of magic labelings of a digraph $D$
as follows. Let $e_{i_1}, \dots, e_{i_{m_i}}$ denote the edges of $D$
that have the vertex $v_i$ as the initial vertex and let $f_{i_1},
f_{i_2}, \dots, f_{i_{s_i}}$ denote the edges of $D$ for which the
vertex $v_i$ is the terminal vertex, then:

\[
\begin{array}{l}
{\cal P}_D = \{ L \in C_D \subseteq \reals^{q}, \hspace{0.05 in}
\sum_{j=1}^{m_i} L(e_{i_j}) \\
 = \sum_{j=1}^{s_i} L(f_{i_j}) = 1 ; i=1,
\dots, n \}.
\end{array}
\]

Let $H_D(r)$ denote the number of magic labelings of $D$ of magic sum
$r$. Then like before, $H_D(r)$ is the Ehrhart quasi-polynomial of
${\cal P}_D$. We now connect the magic labelings of digraphs to magic
labelings of bipartite graphs.

\begin{lemma} \label{digraphbip}
For every digraph $D$, there is a bipartite graph $G_D$
such that the magic labelings of $D$ are in one-to-one correspondence
with the magic labelings of $G_D$. Moreover, the magic sums of the
corresponding magic labelings of $D$ and $G_D$ are also the
same. \end{lemma}

\noindent {\em Proof.}  Denote a directed edge of a digraph $D$ with
$v_i$ as the initial vertex, and $v_j$ as the terminal vertex, by
$e_{ij}$. Let $L$ be a magic labeling of $D$ of magic sum
$r$. Consider a bipartite graph $G_D$ in $2n$ vertices, where the
vertices are partitioned into two sets $A = \{a_1, \dots, a_n \}$ and
$B = \{b_1, \dots, b_n \}$, such that there is an edge between $a_i$
and $b_j$, if and only if, there is an edge $e_{ij}$ in $D$. Consider
a labeling $L_{G_D}$ of $G_D$ such that the edge between the vertices
$a_i$ and $b_j$ is labeled with $L(e_{ij})$.  Observe that the sum of
the labels of the edges incident to $a_{i}$ is the same as the sum of
the labels of incoming edges at the vertex $v_i$ of $D$. Also, the sum
of the labels of edges at a vertex $b_j$ is the sum of the labels of
outgoing edges at the vertex $v_j$ of $D$. Since $L$ is a magic
labeling, it follows that $L_{G_D}$ is a magic labeling of $G_D$ with
magic sum $r$.  Going back-wards, consider a magic labeling
$L^{\prime} $ of $G_D$. We label every edge $e_{ij}$ of $D$ with the
label of the edge between $a_i$ and $b_j$ of $G_D$ to get a magic
labeling $L_D$ of $D$. Observe that $L^{\prime} $ and $L_D$ have the
same magic sum. Hence, there is a one-to-one correspondence between
the magic labelings of $D$ and the magic labelings of $G_D$. $\square$

\begin{figure}
 \begin{center} 
 \includegraphics[scale=0.4]{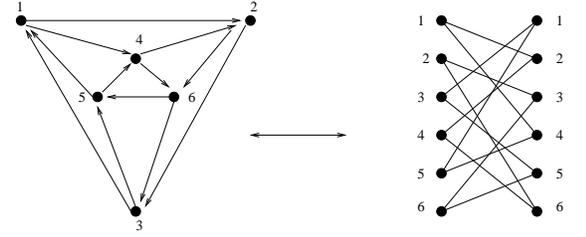}
\caption{A digraph $D_O$ and its corresponding bipartite graph
$G_{D_O}$.} \label{octexample}
\end{center} \end{figure}

For example, the magic labelings of the Octahedral digraph with the
given orientation $D_O$ in Figure \ref{octexample} are in one-to-one
correspondence with the magic labelings of the bipartite graph
$G_{D_O}$.

A graph $G$ is called a {\em positive
graph} if for any edge $e$ of $G$ there is a magic labeling $L$ of $G$
for which $L(e) > 0$ \cite{stanley3}.  Since edges of $G$ that are
always labeled zero for any magic labeling of $G$ may be ignored to
study magic labelings, we will concentrate on positive graphs in
general. We use the following results by Stanley from \cite{stanley3}
and \cite{stanley4} to prove Theorems \ref{graphpolytopethm} and
\ref{digraphpolytopethm} and Corollary \ref{digraphenumerate}.

\begin{thm} [Theorem 1.1, \cite{stanley4}] \label{quasipoly}
Let $G$ be a finite positive graph. Then either $H_G(r)$ is the
Kronecker delta $\delta_{0r}$ or else there exist polynomials $I_G(r)$
and $J_G(r)$ such that $H_G(r) = I_G(r) + (-1)^r J_G(r)$ for all $r
\in \naturals$.
\end{thm}

\begin{thm} [Theorem 1.2, \cite{stanley4}]  \label{degreeofG}
Let $G$ be a finite positive graph with at least one edge. The degree
of $H_G(r)$ is $q-n+b$, where $q$ is the number of edges of $G$, $n$
is the number of vertices, and $b$ is the number of connected
components of $G$ which are bipartite.
\end{thm}

\begin{thm}[Theorem 1.2, \cite{stanley3}] \label{bippoly}
Let $G$ be a finite positive bipartite graph with at least one edge,
then $H_G(r)$ is a polynomial.
\end{thm}

We now conclude that $H_{D}(r)$ is a polynomial for every digraph $D$.

\begin{corollary}  \label{digraphenumerate}
If $D$ is a digraph, then $H_{D}(r)$ is a polynomial of degree
$q-2n+b$, where $q$ is the number of edges of $D$, $n$ is the number
of vertices, and $b$ is the number of connected components of the
bipartite graph $G_D$.
\end{corollary}
\noindent {\em Proof.} The one-to one correspondence between the magic
labelings of $D$ and the magic labelings of $G_D$, implies by Theorem
\ref{bippoly} that $H_D(r)$ is a polynomial, and by Theorem
\ref{degreeofG} that the degree of $H_D(r)$ is $q-2n+b$, where $b$ is
the number of connected components of $G_D$ that are
bipartite. $\square$

Consider the polytope ${\cal P} := \{x|Ax \leq b \}$.  Let $c$ be a
nonzero vector, and let $\delta = $ max $\{cx|Ax \leq b\}$.  The
affine hyperplane $\{x|cx = \delta \}$ is called a {\em supporting
hyperplane} of $\cal P$. A subset $F$ of $\cal P$ is called a {\em
face} of $\cal P$ if $F = \cal P$ or if $F$ is the intersection of
$\cal P$ with a supporting hyperplane of $\cal P$. Alternatively, 
$F$ is a face of $\cal P$ if and only if $F$ is nonempty and 
\[
F = \{ x \in {\cal P} | A^{\prime} x = b^{\prime} \}
\]
for some subsystem $ A^{\prime} x \leq b^{\prime}$ of $Ax \leq b$. See
\cite{schrijver} for basic definitions with regards to polytopes.

Therefore, a  face of ${\cal P}_G$ is a polytope of the form

\[
 \{ L \in {\cal P}_G,   L(e_{i_k}) = 0 ;  e_{i_k} \in E_0 \}, 
\]

where $E_0 = \{ e_{i_1}, \dots, e_{i_r} \}$ is a subset of the set of
edges of $G$.

\begin{thm} \label{graphpolytopethm}
Let $G$ be a finite positive graph with at least one edge. Then the
polytope of magic labelings of $G$, ${\cal P}_G$ is a rational
polytope with dimension $q-n+b$, where $q$ is the number of edges of
$G$, $n$ is the number of vertices, and $b$ is the number of connected
components of $G$ that are bipartite. The $d$-dimensional faces of
${\cal P}_G$ are the $d$-dimensional polytopes of magic labelings of
positive subgraphs of $G$ with $n$ vertices and at most $n-b+d$ edges.
\end{thm}

\begin{figure}
 \begin{center}
     \includegraphics[scale=0.4]{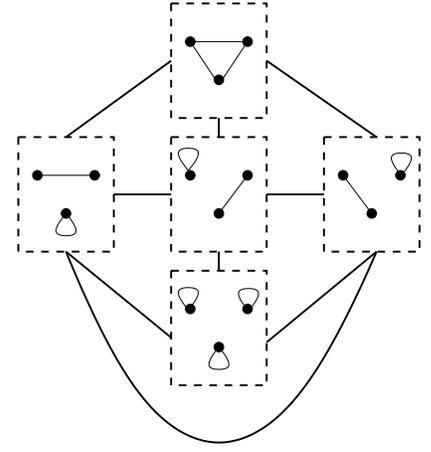}
\caption{The edge graph of  ${\cal P}_{\Gamma_3}$.}  \label{edgegraph}
 \end{center}
 \end{figure}

\begin{figure}
 \begin{center}
     \includegraphics[scale=0.2]{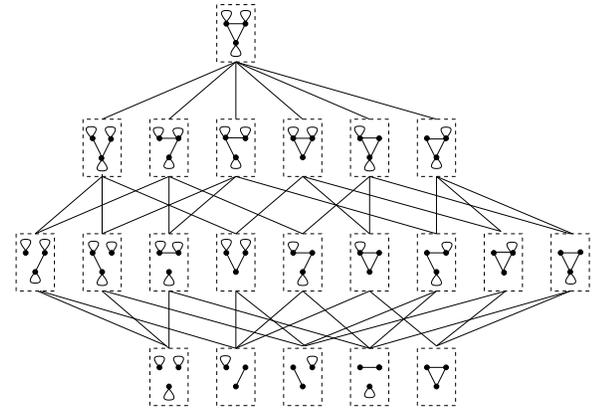}
\caption{The face poset of  ${\cal P}_{\Gamma_3}$.}  \label{faceposet}
 \end{center}
 \end{figure}

Observe from Theorem \ref{graphpolytopethm} that there is an edge
between two vertices $v_i$ and $v_j$ of ${\cal P}_G$ if and only if
there is a graph with at most $n-b+1$ edges, with magic labelings
$v_i$ and $v_j$.  The edge graph of ${\cal P}_{\Gamma_3}$ is given in
Figure \ref{edgegraph}. Similarly, we can draw the face poset of
${\cal P}_{G}$ (see Figure \ref{faceposet} for the face poset of
${\cal P}_{\Gamma_3}$).

An $n \times n$ {\em semi-magic square} of magic sum $r$ is an $n
\times n$ matrix with nonnegative integer entries such that the
entries of every row and column add to $r$. {\em Doubly stochastic
matrices} are $n \times n$ matrices in $\reals^{n^2}$ such that their
rows and columns add to 1. The set of all $n \times n$ doubly
stochastic matrices form a polytope $B_n$, called the {\em Birkhoff
polytope}. See \cite{billerasarang}, \cite{brualdigibson}, or
\cite{schrijver} for a detailed study of the Birkhoff polytope.

A {\em symmetric magic square} is a semi-magic square that is also a
 symmetric matrix. Let $H_n(r)$ denote the number of symmetric magic
 squares of magic sum $r$ (see \cite{mayathesis}, \cite{stanley4}, and
 the references therein for the enumeration of symmetric magic
 squares).  We define the polytope ${\cal S}_n$ of $n \times n$
 symmetric magic squares to be the convex hull of all real nonnegative
 $n \times n$ symmetric matrices such that the entries of each row
 (and therefore column) add to one.

\begin{figure}
 \begin{center}
     \includegraphics[scale=0.5]{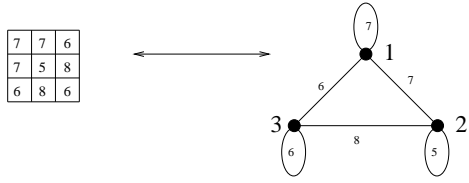}
\caption{A magic labeling of $\Gamma_3$ and its corresponding
symmetric magic square.} \label{gamma3eg}
 \end{center}
 \end{figure}

 A one-to-one correspondence between symmetric magic squares $M =
 [m_{ij}]$ of magic sum $r$, and magic labelings of the graph
 $\Gamma_n$ of the same magic sum $r$ was established in
 \cite{stanley4}: let $e_{ij}$ denote an edge between the vertex $v_i$
 and the vertex $v_j$ of $\Gamma_n$. Label the edge $e_{ij}$ of
 $\Gamma_n$ with $m_{ij}$, then this labeling is a magic labeling of
 $\Gamma_n$ with magic sum $r$. See Figure \ref{gamma3eg} for an
 example. Therefore, we get ${\cal P}_{\Gamma_n}$ = ${\cal S}_{n}$ and
 $H_{\Gamma_n}(r) = H_n(r)$.

\begin{corollary} \label{mayathm}
The polytope of magic labelings of the complete general graph ${\cal
P}_{\Gamma_n}$ is an $n(n-1)/2$ dimensional rational polytope with the
following description
\[{\cal P}_{\Gamma_n} = \left  \{ \begin{array}{l}
 L = (L(e_{ij}) \in {\reals}^{\frac{n(n+1)}{2}}; \\
  L(e_{ij}) \geq 0; \\
1 \leq i,j \leq n, i \leq j, \\     
\sum_{j=1}^i L(e_{ji}) + \sum_{j=i+1}^n L(e_{ij}) = 1 \\
\mbox{ for } i=1, \dots, n 
\end{array}
\right \}
\]
The $d$-dimensional faces of ${\cal P}_{\Gamma_n}$ are $d$-dimensional
polytopes of magic labelings of positive graphs with $n$ vertices and
at most $n+d$ edges. There are $2n-1 \choose n$ faces of ${\cal
P}_{\Gamma_{2n}}$ that are copies of the Birkhoff polytope $B_n$.
\end{corollary}

We define a digraph $D$ to be a {\em positive digraph} if the corresponding 
bipartite graph $G_D$ is positive. 

\begin{thm} \label{digraphpolytopethm}
Let $D$ be a positive digraph with at least one edge. Then, ${\cal
P}_D$ is an integral polytope with dimension $q-2n+b$, where $q$ is
the number of edges of $D$, $n$ is the number of vertices, and $b$ is
the number of connected components of $G_D$ that are bipartite. The
$d$-dimensional faces of ${\cal P}_D$ are the $d$-dimensional
polytopes of magic labelings of positive subdigraphs of $D$ with $n$
vertices and at most $2n-b+d$ edges.
\end{thm}

\begin{figure}
 \begin{center} 
 \includegraphics[scale=0.3]{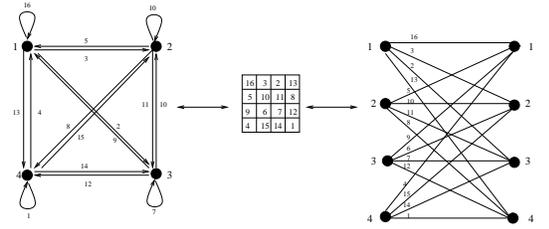}
\caption{Two different graph labelings associated to a semi-magic square.}
 \label{semi-magic} \end{center} \end{figure}

Let $\Pi_n$ denote the complete digraph with $n$ vertices, i.e, there
is an edge from each vertex to every other, including the vertex
itself (thereby creating a loop at every vertex), then $G_{\Pi_n}$ is
the the complete bipartite graph $K_{n,n}$. We get a one-to-one
correspondence between semi-magic squares $M = [m_{ij}]$ of magic sum
$r$ and magic labelings of $\Pi_n$ of the same magic sum $r$ by
labeling the edges $e_{ij}$ of $\Pi_n$ with $m_{ij}$.  This also
implies that there is a one-to-one correspondence between semi-magic
squares and magic labelings of $K_{n,n}$ (this correspondence is also
mentioned in \cite{stanley3} and \cite{stewart}).  See Figure
\ref{semi-magic} for an example.

A good description of the faces of Birkhoff polytope is not known \cite{pak}.
We can now give an explicit description of the faces of the Birkhoff polytope.

\begin{corollary} \label{facesbirk}
${\cal P}_{\Pi_n}$ is the Birkhoff polytope $B_n$.  The
$d$-dimensional faces of $B_n$ are polytopes of magic labelings of
positive digraphs with dimension $d$, $n$ vertices and at most
$2n+d-1$ edges. The vertices of ${\cal P}_D$, where $D$ is a positive
digraph, are permutation matrices.
 
\end{corollary}

\begin{figure}
 \begin{center}
     \includegraphics[scale=0.4]{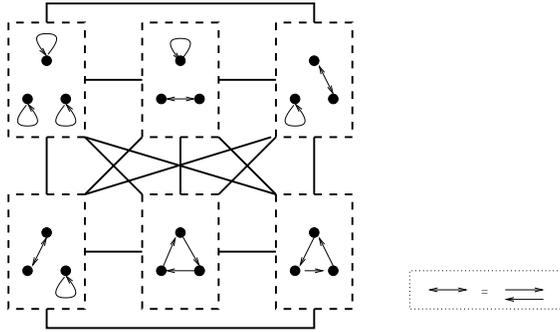}
\caption{The edge graph of the Birkhoff Polytope $B_3$.}  \label{b3edgegraph}
 \end{center}
 \end{figure}

\begin{figure}
 \begin{center}
     \includegraphics[scale=0.4]{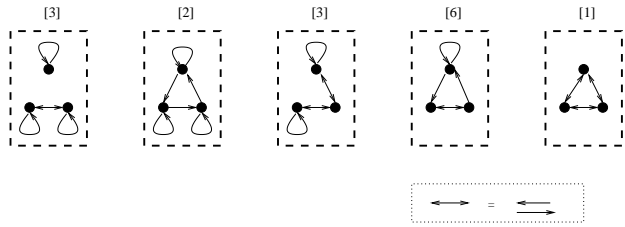}
\caption{The generators of the edges of the Birkhoff
Polytope $B_3$.}  \label{b3edges}
 \end{center}
 \end{figure}

\begin{figure}
 \begin{center}
     \includegraphics[scale=0.4]{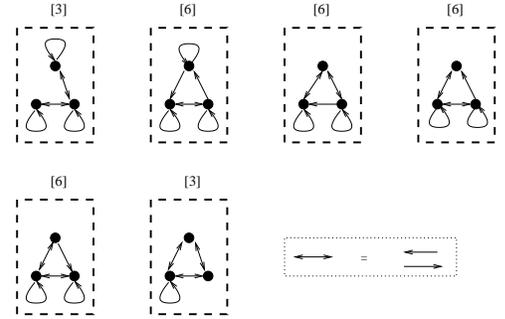}
\caption{The generators of the 2-dimensional faces of the Birkhoff
Polytope $B_3$.}  \label{b32dimfaces}
 \end{center}
 \end{figure}

\begin{figure}
 \begin{center}
     \includegraphics[scale=0.4]{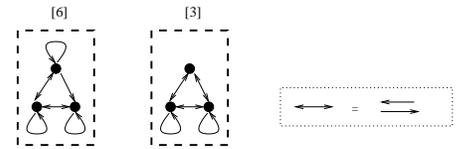}
\caption{The generators of the facets of the Birkhoff
Polytope $B_3$.}  \label{b3facets}
 \end{center}
 \end{figure}

\begin{figure}
 \begin{center}
     \includegraphics[scale=0.4]{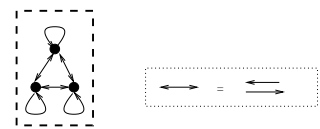}
\caption{The Birkhoff Polytope $B_3$.}  \label{b3}
 \end{center}
 \end{figure}

 See Figure \ref{b3edgegraph} for the edge graph of $B_3$. Two faces
 of a polytope of magic labelings of a graph (or a digraph) are said
 to be {\em isomorphic faces} if the subgraphs (subdigraphs,
 respectively) defining the faces are isomorphic. A set of faces is
 said to be a {\em generating set of $d$-dimensional faces} if every
 $d$-dimensional face is isomorphic to one of the faces in the
 set. See Figures \ref{b3edges}, \ref{b32dimfaces}, \ref{b3facets},
 and \ref{b3} for the generators of the edges, the two dimensional
 faces, the facets, and the Birkhoff polytope $B_3$, respectively.
 The numbers in square brackets in the figures indicate the number of
 faces in the isomorphism class of the given face. 

The proofs of Theorems \ref{graphpolytopethm} and
\ref{digraphpolytopethm}, and Corollaries \ref{mayathm} and
\ref{facesbirk} are presented in Section \ref{polytopesection}.

\section{Polytopes of magic labelings.} \label{polytopesection}
A polytope {$\cal P$} is called {\em rational} if each vertex of
{$\cal P$} has rational coordinates. 

An element $\beta$ in the semigroup $S_{C_G}$ is said to be {\em
completely fundamental}, if for any positive integer $n$ and $\alpha,
{\alpha}^{\prime} \in S_{C_G}$, $n \beta = \alpha + {\alpha}^{\prime}$
implies $\alpha = i \beta$ and ${\alpha}^{\prime} = (n-i) \beta$, for
some positive integer $i$, such that $0 \leq i \leq n$ (see
\cite{stanley2}).

\begin{lemma} \label{polytopeverticeslemma} ${\cal P}_G$ is a rational 
polytope.
\end{lemma}
\noindent {\em Proof.} Proposition $4.6.10$ of Chapter 4 in
\cite{stanley2} states that the set of extreme rays of a cone and the
set of completely fundamental solutions are identical. Proposition 2.7
in \cite{stanley3} states that every completely fundamental magic
labeling of a graph $G$ has magic sum 1 or 2. Thus, the extreme rays
of the cone of magic labelings of a graph $G$ are irreducible
2-matchings of $G$. We get a vertex of ${\cal P}_{G}$ by dividing the
entries of a extreme ray by its magic sum. Thus, ${\cal P}_G$ is a
rational polytope.  $\square$

\begin{lemma} \label{polytopedimensionlemma}
The dimension of ${\cal P}_G$ is $q-n+b$, where $q$ is the number of
edges of $G$, $n$ is the number of vertices, and $b$ is the number of
connected components that are bipartite.
\end{lemma}

\noindent {\em Proof.}  Ehrhart's theorem states that the degree of
$H_G(r)$ is the dimension of ${\cal P}_G$ \cite{becketal}. The degree
of $H_G(r)$ is $q-n+b$ by Theorem \ref{degreeofG}. Therefore, the
dimension of ${\cal P}_G$ is $q-n+b$.  $\square$

\begin{lemma} \label{polytopefaceslemma}
The $d$-dimensional faces of ${\cal P}_G$ are the $d$-dimensional
polytopes of magic labelings of positive subgraphs of $G$ with $n$
vertices and at most $n-b+d$ edges.
\end{lemma}
\noindent {\em Proof.} An edge $e$ labeled with a zero in a magic
labeling $L$ of $G$ does not contribute to the magic sum, therefore,
we can consider $L$ as a magic labeling of a subgraph of $G$ with the
edge $e $ deleted. Since a face of ${\cal P}_{G}$ is the set of magic
labelings of $G$ where some edges are always labeled zero, it follows
that the face is also the set of all the magic labelings of a subgraph
of $G$ with these edges deleted.  Similarly, every magic labeling of a
subgraph $H$ with $n$ vertices corresponds to a magic labeling of $G$,
where the missing edges of $G$ in $H$ are labeled with 0. Now, let $H$
be a subgraph such that the edges $e_{r1}, \dots, e_{rm}$ are labeled
zero for every magic labeling of $H$. Then the face defined by $H$ is
same as the face defined by the positive graph we get from $H$ after
deleting the edges $e_{r1}, \dots, e_{rm}$. Therefore, the faces of
${\cal P}_G$ are polytopes of magic labelings of positive subgraphs.

By Lemma \ref{polytopedimensionlemma}, the dimension of ${\cal P}_G$
is $q-n+b$. Therefore, to get a $d$-dimensional polytope, we need to
label at least $q-n+b-d$ of $G$ edges always 0. This implies that the
$d$-dimensional face is the set of magic labelings of a positive
subgraph of $G$ with $n$ vertices and at most $n-b+d$ edges.
$\square$

The proof of Theorem \ref{graphpolytopethm} follows from Lemmas
\ref{polytopeverticeslemma}, \ref{polytopedimensionlemma}, and
\ref{polytopefaceslemma}. We can now prove Corollary \ref{mayathm}.

\noindent {\em Proof of Corollary \ref{mayathm}.}  

It is clear from the one-to-one correspondence between magic labelings
of $\Gamma_n$ and symmetric magic squares that ${\cal P}_{\Gamma_n}$
has the given description. Since the graph ${\Gamma_n}$ has
$\frac{n(n+1)}{2}$ edges and $n$ vertices, and every graph is a
subgraph of $\Gamma_n$, it follows from Theorem \ref{graphpolytopethm}
that the dimension of ${\cal P}_{\Gamma_n}$ is $\frac{n(n-1)}{2}$; the
$d$-dimensional faces of ${\cal P}_{\Gamma_n}$ are $d$-dimensional
polytopes of magic labelings of positive graphs with $n$ vertices and
at most $n+d$ edges.

We can partition the vertices of $\Gamma_{2n}$ into two equal sets $A$
and $B$ in ${2n-1} \choose n$ ways: Fix the vertex $v_1$ to be in the
set $A$ , then we can choose the n vertices for the set $B$ in ${2n-1}
\choose n$ ways, and the remaining $n-1$ vertices will belong to the
set $A$. By adding the required edges, we get a complete bipartite
graph for every such partition of the vertices of $\Gamma_{2n}$. Thus,
the number of subgraphs of $\Gamma_{2n}$ that are isomorphic to
$K_{n,n}$ is ${2n-1} \choose n$. Therefore, there are $2n-1 \choose n$
faces of ${\cal P}_{\Gamma_{2n}}$ that are Birkhoff polytopes because
every isomorphic copy of $K_{n,n}$ contributes to a face of ${\cal
P}_{\Gamma_{2n}}$.  $\square$

\begin{lemma}
Let $G$ be a graph with $n$ vertices. A labeling $L$ of $G$ with magic
sum $s$ can be lifted to a magic labeling $L^{\prime}$ of the complete
general graph $\Gamma_n$ with magic sum $s$.
\end{lemma}

\noindent {\em Proof.}
Since $G$ is a subgraph of $\Gamma_n$, every
labeling $L$ of $G$ can be lifted to a labeling $L^{\prime}$ of
$\Gamma_n$, where
\[
L^{\prime} (e_{ij}) = \left \{ \begin{array}{ll} L(e_{ij}) & \mbox{ if
$e_{ij}$ is also an edge of $G$,} \\ 0 & \mbox{otherwise.}
\end{array}
\right .
\]

Since  the edges with nonzero labels are the same for both $L$  and 
$L^{\prime}$, it follows that the magic sums are also the same. $\square$ 

\begin{lemma} \label{hilbGprop}
Let $G$ be a graph with $n$ vertices. The minimal Hilbert basis of
$C_G$ can be lifted to a subset of the minimal Hilbert basis of
$C_{\Gamma_n}$.
\end{lemma}
\noindent {\em Proof.}  If $L$ is an irreducible magic labeling of
$G$, then clearly it lifts to an irreducible magic labeling $L^{\prime}$ of
$\Gamma_n$. Since the minimal Hilbert basis is the set of all
irreducible magic labelings, we get that the minimal Hilbert basis of
$C_G$ corresponds to a subset of the minimal Hilbert basis of
$C_{\Gamma_n}$.  $\square$

Similarly, we can prove:
\begin{lemma} \label{digraphhilb} 
For a digraph $D$ with $n$ vertices, a magic labeling $L$ with magic
sum $s$ can be lifted to a magic labeling $L^{\prime}$ of $\Pi_n$ with
the same magic sum $s$. The minimal Hilbert basis of $C_D$ can be
lifted to a subset of the minimal Hilbert basis of $C_{\Pi_n}$.
\end{lemma}

\begin{lemma} \label{magicsum1}
Let $D$ be a digraph with $n$ vertices. All the elements of the
minimal Hilbert basis of $C_D$ have magic sum 1.
\end{lemma}
\noindent {\em Proof.}  It is well-known that the minimal Hilbert
basis of semi-magic squares are the permutation matrices (see
\cite{schrijver}) and therefore have magic sum 1. The one-to-one
correspondence between magic labelings of $\Pi_n$ and semi-magic
squares implies that the minimal Hilbert basis elements of $C_{\Pi_n}$
have magic sum 1. It follows by Lemma \ref{digraphhilb} that all the
elements of the minimal Hilbert basis of $C_{D}$ have magic sum
1. $\square$

We now prove our results about polytope of magic digraphs.

\noindent {\em Proof of Theorem \ref{digraphpolytopethm}.}  By Lemma
\ref{magicsum1}, all the elements of the Hilbert basis of $C_D$ have
magic sum 1. Since the extreme rays are a subset of the Hilbert basis
elements, it follows that the vertices of ${\cal P}_D$ are
integral. Since ${\cal P}_D = {\cal P}_{G_D}$, it follows by Theorem
\ref{graphpolytopethm} that the dimension of ${\cal P}_D$ is $q-2n+b$;
the $d$-dimensional faces of ${\cal P}_{D}$ are the $d$-dimensional
polytopes of magic labelings of positive subdigraphs of $D$ with $n$
vertices and at most $2n-b+d$ edges.  $\square$

We derive our results about the faces of the Birkhoff polytope as a
consequence.

{\em Proof of Corollary \ref{facesbirk}.}  The one-to-one
correspondence between semi-magic squares and magic labelings of
${\Pi_n}$ gives us that ${\cal P}_{\Pi_n} = B_n$. Since every digraph
with $n$ vertices is a subgraph of $\Pi_n$, by Theorem
\ref{digraphpolytopethm}, it follows that its $d$-dimensional faces
are $d$-dimensional polytopes of magic labelings of positive digraphs
with $n$ vertices and at most $2n-1+d$ edges. Since the vertex set of
a face of $B_n$ is a subset of the vertex set of $B_n$ it follows that
the vertices of ${\cal P}_D$, where $D$ is a positive digraph, are
permutation matrices. $\square$

\section{Applications.}
\begin{figure}
 \begin{center}
\includegraphics[scale=0.5]{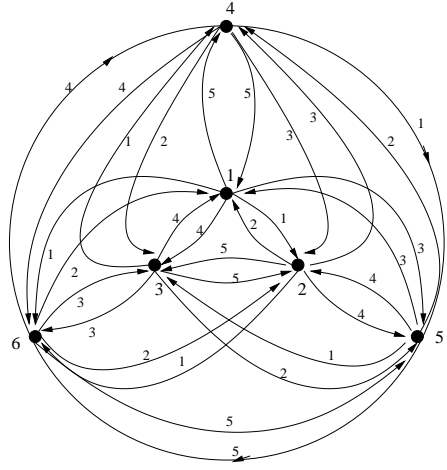}
\caption{Cayley digraph of the group $S_3$ \cite{konig}.}
 \label{cayleys3} \end{center} \end{figure}

In this section, we present some examples and applications of magic
graphs.  Interesting examples of magic digraphs are Cayley digraphs of
finite groups. Let $G$ be a finite group $\{g_1,g_2, \dots, g_n= I
\}$. The {\em Cayley group digraph} of $G$ is a graphical
representation of $G$: every element $g_i$ of the group $G$
corresponds to a vertex $v_i$ $(i=1,2, \dots, n)$ and every pair of
distinct vertices $v_i, v_j$ is joined by an edge labeled with
$\alpha$ where $g_{\alpha} = g_j g_i^{-1}$ \cite{konig}. For example,
the Cayley digraph for the permutation group
\[
\begin{array}{l}
S_3 = \{ g_1 = (123), g_2 = (132), g_3 = (23), \\
 g_4 = (12), g_5 = (13), g_6 = I \}
\end{array}
\]
 is given in Figure \ref{cayleys3}. 

\begin{prop} \label{cayley}
The Cayley digraph of a group of order $n$ is a magic digraph with
magic sum $\frac{n(n-1)}{2}$.
\end{prop}

\noindent {\em Proof.} Let $e_{ij}$ denote an edge between the vertex
$v_i$ and $v_j$ of the Cayley digraph such that $v_i$ is the initial
vertex and $v_j$ is the terminal vertex.  Let $v_l$ be a vertex of the
Cayley digraph, and let $\alpha$ be an integer in the set $\{ 1, 2,
\dots , n-1 \}$.  Let $g_p = g_{\alpha}g_l$ and let $g_q = g_l
g_{\alpha}$. Then, the edges $e_{lp}$ and $e_{ql}$ are labeled by
$\alpha$. Also, $g_jg_i^{-1} = g_n = I$ if and only if $i=j$.  Hence,
a Cayley group digraph is a magic digraph with magic sum $1 + 2 +
\cdots + (n-1) = \frac{n(n-1)}{2}$ (see also Chapter 8, Section 5 in
\cite{konig}). $\square$

A digraph is called {\em Eulerian} if for each vertex $v$ the indegree
and the outdegree of $v$ is the same. Therefore, Eulerian digraphs can
also be studied as magic digraphs where all the edges are labeled by 1
(see \cite{alantarsi} for the applications of Eulerian digraphs to
digraph colorings).
An {\em $n$-matching} of $G$ is a magic labeling of $G$ with magic sum
at most $n$ and the labels are from the set $\{0, 1,\dots,n \}$ (see
\cite{lovasz}, chapter 6).  A {\em perfect matching} of $G$ is a
1-matching of $G$ with magic sum 1.

\begin{prop} \label{perfectmatchings} 
The perfect matchings of $G$ are the minimal Hilbert basis elements of
$C_G$ of magic sum 1 and the number of perfect matchings of $G$ is
$H_G(1)$. 
\end{prop}

\noindent{\em Proof.}  Magic labelings of magic sum 1 always belong to
the minimal Hilbert basis because they are irreducible.  Therefore,
perfect matchings belong to the minimal Hilbert basis because they
have magic sum 1. Conversely, every magic labeling of magic sum 1 is a
perfect matching. So we conclude that the perfect matchings of $G$ are
the minimal Hilbert basis elements of $C_G$ of magic sum 1. The fact
that the number of perfect matchings of $G$ is $H_G(1)$ follows by the
definition of $H_G(1)$.  $\square$

\begin{figure}
 \begin{center} 
 \includegraphics[scale=0.5]{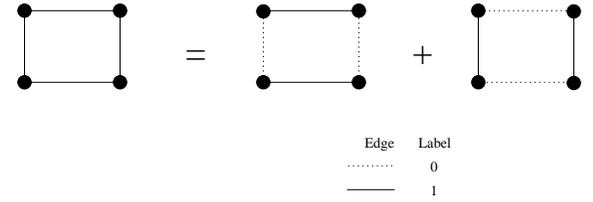}
\caption{Graph Factorization.} \label{factor}
\end{center}
\end{figure}

Hilbert basis can also be used to study factorizations of labeled
graphs. We define {\em Factors} of a graph $G$ with a labeling $L$ to
be labelings $L_i, i = 1, \dots, r$ of $G$ such that $L(G) =
\sum_{i=1}^r L_i(G)$, and if $L_i(e_k) \neq 0$ for some edge $e_k$ of
$G$, then $L_j (e_k) = 0$ for all $ j \neq i$. A decomposition of $L$
into factors is called a {\em factorization} of $G$.  An example of a
graph factorization is given in Figure \ref{factor}. See Chapters 11
and 12 of \cite{konig} for a detailed study of graph factorizations.


Our results enable us to reprove some known facts about the Birkhoff
polytope as well.  For example, Theorem \ref{digraphpolytopethm} gives
us that the dimension of $B_n$ is ${(n-1)}^2$. The leading coefficient
of the Ehrhart polynomial of $B_n$ is the volume of $B_n$. This number
has been computed for $n=1,2, \dots, 9$ (see \cite{beck} and
\cite{chanrobbins}). 

The software 4ti2 \cite{raymond} can be used to find the Hilbert bases of
$C_G$, and the software LattE \cite{rudy} can be used to compute the
generating functions of $H_G(r)$ effectively.  See \cite{mayathesis}
for other results about magic labelings of graphs. See \cite{wallis}
for a study of other types of magic graphs.
   
\noindent {\bf Acknowledgments.} I thank my advisor Jes\'us De Loera
for several suggestions and useful conversations.

\bibliographystyle{plain}

\end{document}